\documentclass[12pt]{amsart}

\usepackage[utf8]{inputenc}

\usepackage[T1]{fontenc}

\usepackage[a4paper, margin=2.7cm]{geometry}

\usepackage{amsmath, amsthm, amsfonts, amssymb}
\usepackage{mathrsfs}           

\usepackage[colorlinks=true,linkcolor=blue,citecolor=blue,urlcolor=blue,breaklinks]{hyperref}

\usepackage{bbm}

 \newcommand{\R}{{\mathbb R}}

\def\Om{\Omega}

\usepackage{xcolor}

\theoremstyle{definition}

\theoremstyle{remark}

\theoremstyle{example}

\def\R{\mathbb R}

\title{Integro-differential equations with delays: A perturbation approach}
\author{Hamid Bounit}
\address{Department of Mathematics, Faculty of Sciences Hay Dakhla, BP8106, 80000--Agadir, Morocco}
\email{h.bounit@uiz.ac.ma}
\author{Abderrahim Driouich}
\address{Department of Mathematics, Faculty of Sciences Hay Dakhla, BP8106, 80000--Agadir, Morocco}
\email{a.driouich@uiz.ac.ma}
\author{Said Hadd}
\address{Department of Mathematics, Faculty of Sciences Hay Dakhla, BP8106, 80000--Agadir, Morocco}
\email{s.hadd@uiz.ac.ma}

\begin{document}
\maketitle
\begin{abstract}
This paper focuses on the study of integro-differential equations with delays, presenting a novel perturbation approach. The primary objective is to introduce the concepts of classical and mild solutions for these equations and establish their existence and uniqueness, under suitable assumptions. Furthermore, we provide a variation of constant formula that characterizes these solutions. To illustrate the applicability of the proposed methodology, we present an example of integro-differential Volterra equations with a nonlocal kernel. In addition to the aforementioned contributions, a secondary goal of this paper is to address an issue concerning the statement and proof of a fundamental theorem presented in a previous work \cite{Zaza}. Specifically, we aim to rectify the statement and provide a corrected proof for Theorem 2.6 in \cite{Zaza}. By doing so, we enhance the accuracy and reliability of the existing literature in this field.
\end{abstract}
\section{Introduction}\label{sec:1}
Integro-differential Volterra equations with delays are a type of functional differential equation that combines integral and differential terms with delay effects. They are named after Vito Volterra \cite{Volterra}, an Italian mathematician who made significant contributions to the field of mathematical analysis in the late 19th and early 20th centuries. Note that Volterra equations with delays arise in many areas of science and engineering, including physics, biology, economics, and control theory.

This paper focuses on establishing the existence, uniqueness, and variation of constants formula for both mild and classical solutions of the following integro-differential equation with delays
\begin{align}\label{Volterra-f}
\begin{split}
& \dot{u}(t)=Au(t)+\int^t_0 a(t-s)Cu(s)ds\cr &\hspace{3cm}+\int^t_0 b(t-s)Lu_sds+Ku_t+f(t),\quad t\ge 0,\cr & u(0)=x,\quad u(t)=\varphi(t),\quad t\in [-1,0],
\end{split}
\end{align}
Here $A:D(A)\subset X\to X$ is the generator of a strongly continuous semigroup $T:=(T(t))_{t\ge 0}$ on a Banach space $X$, $C:D(A)\to X$ is a linear (not closed nor closable) operator on $X$, the delay operator $L,K:W^{1,p}([-1,0],X)\to X$ are linear,  $a(\cdot)$ and $b(\cdot)$ are scaler functions belonging to the Lebesgue space $L^{p}(\R^+)$, the nonhomogeneous term $f\in L^p(\R^+,X)$ for some $p\in (1,\infty),$ and the initial conditions $x\in X$ and $\varphi\in L^p([-1,0],X)$.
Here $u(\cdot):[-1,\infty)\to X$ is the eventual solution to \eqref{Volterra-f}. At any time $t\ge 0,$ the history function of $u(\cdot)$  is the function $u_t:[-1,0]\to X$ with $u_t(\theta)=u(t+\theta)$ for any $\theta\in [-1,0]$.

This paper introduces a perturbation approach to the equation \eqref{Volterra-f}, presenting a novel method for analyzing and solving these equations. By incorporating perturbation techniques, we aim to derive a new variation of constant formula for the solutions, utilizing Yosida extensions of the  operators $K,L$, and $C$. The perturbation approach offers several advantages over classical methods for analyzing integro-differential equations with delays. Firstly, by considering perturbations, we gain insights into the effects of small changes or disturbances on the solutions of the equation \eqref{Volterra-f}. This allows us to understand the system's response under perturbations and quantify their impact using the derived variation of constants formula. Secondly, by utilizing Yosida extensions, we introduce a powerful mathematical tool that helps us handle the delay operators $K,L$ and the kernel operator $C$ in a unified framework. This enables a more comprehensive analysis and solution of the integro-differential equations with delays.

Free-delay integro-differential equations have been extensively studied in the literature due to their wide applicability and significance in various fields. In this paper, we highlight key references that are relevant to the investigation of these equations, namely \cite{DS}, \cite[Chap.VI, Section 7(c)]{EngNag}, and \cite{Pruss}. These references provide valuable insights and utilize a semigroup approach to solve and analyze free-delay integro-differential equations. The reference \cite{DS} offers significant contributions to the study of free-delay integro-differential equations. It presents in-depth analysis and provides fundamental results related to the existence, uniqueness, and stability of solutions. The chapter \cite[Chap.VI, Section 7(c)]{EngNag} offers a comprehensive overview of the theory, highlighting the challenges and intricacies associated with free-delay integro-differential equations in the case of $C=A$ and $L=K=0$. The authors provide rigorous mathematical formulations and techniques, enabling a deeper understanding of these equations and their solutions.

The book \cite{Brun-02} focuses on collocation methods for solving integro-differential equations, including those with delays. It provides a thorough treatment of numerical methods and their convergence analysis, along with practical examples and applications. The reference \cite{LMN-17} explores fractional integro-differential equations with delays, which generalize the classical theory of integro-differential equations. It covers existence, uniqueness, stability, and numerical methods for fractional integro-differential equations with delays. Although focused on delay differential equations, the reference \cite{GL-91} also covers some aspects of integro-differential equations with delays. It provides a detailed analysis of oscillatory behavior and stability properties of solutions to delay equations, including integral terms. More recently, in the paper \cite{EHB}, the authors introduced an insightful analytic approach to address the equation \eqref{Volterra-f} specifically for the scenario where $C=A$ and $L=0$. They successfully derived a variation of constants formula by employing the concept of resolvent families linked to the free-delay integro-differential equation. Additionally, their work incorporates a comprehensive and detailed spectral theory, along with a significant contribution to control theory.

In recent years, significant advancements have been made in the study of integro-differential equations with delays, specifically focusing on Equation \eqref{Volterra-f}. One notable contribution is the work by the authors of \cite{EHB}, where they presented an analytical approach to this equation in the case of $C=A$ and $L=0$. Their study encompasses a variation of constant formula derived through the utilization of resolvent families associated with the free-delay integro-differential equation. Additionally, they introduce a detailed spectral theory and explore the application of control theory.

\textcolor{blue}{More recently, a work \cite{Zaza} has emerged, addressing Equation \eqref{Volterra-f} under the conditions of $K=0$ and $f=0$. This work draws inspiration from previous studies such as \cite{ABDH}, \cite{EHB}, and incorporates perturbation techniques using admissible observation operators developed in \cite{Hadd}, which generalize the Miyadera-Voigt perturbation theorem. However, it is important to note that the work \cite{Zaza} contains several issues that we will outline in this discussion.} \textcolor{red} { Prior to undertaking the present work, we took the initiative to contact the journal in which the work \cite{Zaza} was published. Our objective was to highlight the lack of proper citation of our ideas within \cite{Zaza}, as well as to bring attention to several issues present in his study. Unfortunately, the feedback we received from the editorial board of the journal was negative. Consequently, we made the decision to embark on a new research endeavor that aims to rectify the issues of \cite{Zaza} and incorporate recent, intriguing developments in the field.}

 \textcolor{blue} {One issue with the work of \cite{Zaza} is that the statement of \cite[Theorem 2.6]{Zaza} is incomplete and lacks some essential details. Specifically, the condition $x\in D(A)$ needs to be replaced with $x\in D(A)$ and $\varphi\in W^{1,p}([-1,0],X)$ such that $x=\varphi(0)$. Furthermore, the well-posedness of integro-differential equation with delay is not clearly defined before this theorem, even though it is the focus of the theorem. Clarifying the definition of classical solutions before stating} \textcolor{red} {\cite[Theorem 2.6]{Zaza}} \textcolor{blue}{ would improve the clarity of the work  \cite{Zaza}. The variation of constants formula stated in this theorem involves the terms $Cu(\tau)$ and $Lu_\tau$ (the notation in \cite{Zaza} was $C=F$). However, it is unclear why these terms are well-defined, given that the operators $C$ and $L$ are defined on $D(A)$ and $W^{1,p}([-1,0],X)$, respectively. This highlights the importance of first clearly defining the concept of solutions for delayed integro-differential equation before using it to prove the result in}\textcolor{red}{ \cite[Theorem 2.6]{Zaza}}.

\textcolor{red}{The proof of \cite[Theorem 2.6]{Zaza} presents a significant gap.} \textcolor{blue}{Specifically, the proof does not provide a complete explanation of why the classical solution $u(\cdot)$ of \eqref{Volterra-f} exists in the case of $K=0$ and $f=0$. The fourth and fifth pages of \cite{Zaza} contain formulas that are presented in a formal manner, without proper justification. While these formulas provide an idea about the free-delay Cauchy problem and its solution, they do not offer a complete proof of the existence and uniqueness of classical solutions. The author mention that he follow the strategy of the well-known reference \cite{DS}, but it is important to note that in that reference, the authors explicitly state that they use formal equations. Therefore, to improve the rigor of the proof in \cite[Theorem 2.6]{Zaza}, it is necessary to provide a clear and complete justification of the formulas used, and to establish the existence and uniqueness of classical solutions in a rigorous manner. Once the existence of the classical solution is established, the next step is to prove the variation of constants formula for such a solution. While this may seem like a straightforward task, it is important to ensure that the proof is rigorous and clearly explains the steps involved. By carefully proving the variation of constants formula, we can ensure that it is valid for all classical solutions to delay integro-differential equation, and use it to derive useful information about the behavior of the solutions. Therefore, it is essential to provide a well-justified and detailed proof of the variation of constants formula for classical solutions.}

\textcolor{blue}{In addition to the concerns raised regarding the work \cite{Zaza}, it is worth noting that a related study, namely \cite{Zaza1}, also exhibits similar issues and shortcomings. These works share notable similarities in their methodology, results, and deficiencies. Therefore, it is imperative to address and rectify the concerns associated with both papers in order to ensure the reliability and accuracy of the scientific contributions.}

The organization of the present work is as follows: In Section 2, we delve into the topic of admissible observation operators and their connection to perturbed Cauchy problems. Specifically, we provide a comprehensive overview of relevant facts in this area. We recall a perturbation theorem, originally presented in \cite{Hadd}, which serves as a key foundation for our subsequent analysis. Additionally, we present a significant result that demonstrates the invariance of admissibility for observation operators under a specific class of perturbations. These findings contribute to a deeper understanding of the relationship between observation operators and perturbed systems, shedding light on their interplay and implications in the broader context of integro-differential equations with delays.

In Section 3, we provide a comprehensive analysis of the integro-differential equation with delay \eqref{Volterra-f} by introducing and exploring two important solution concepts: classical solutions and mild solutions. Firstly, we define classical solutions to \eqref{Volterra-f} in Definition 3.1 and establish their existence and uniqueness through rigorous proof. This sets the foundation for understanding the behavior and properties of the equation. Next, we introduce the concept of mild solutions in Definition 3.6 and demonstrate that every classical solution also qualifies as a mild solution to \eqref{Volterra-f}. This result bridges the gap between the two solution frameworks and highlights their inherent connections. To further advance our understanding, we present a theorem that outlines the conditions under which a unique mild solution exists for \eqref{Volterra-f}. This result provides valuable insights into the existence and uniqueness aspects of the equation, facilitating a more comprehensive analysis of its behavior. Lastly, we conclude the section by considering an application of these concepts to an integro-differential equation governed by the Laplcian operator on $L^2(\Om)$ for an open bounded set $\Om\subset \R^n$ with smooth boundary $\partial\Om$ subject to Neumann boundary conditions and a nonlocal perturbation kernel $C$. This practical example illustrates the applicability and relevance of our theoretical findings in real-world scenarios, further enhancing the practical significance of our research.

\newpage

{\Huge The rest of the content of the article will be given here once published in a journal. We are forced not to put all the details of the results at this stage because we have already had bad experiences lately, as already explained... But certainly, in a few weeks, we will put the complete paper in Arxiv. }

\end{document}